\documentclass[a4paper,12pt,leqno]{article}
\usepackage{latexsym}
\usepackage[all]{xy}

\usepackage{amssymb} % \mathbb
\usepackage{amsmath} % \mathbb
\usepackage{theorem}

\usepackage{theorem}
\usepackage{amscd}
\usepackage{latexsym}
\usepackage{enumerate}
\usepackage[dvips]{graphicx,color}

\def\R{{\mathbb{R}}}% \R == \mathbb{R}
\def\A{{\mathcal{A}}}% \A == \mathcal{A}
\def\Shi{{\mathcal{S}}}% \B == \mathcal{B}
\numberwithin{equation}{section}

\newtheorem{theorem}{Theorem}[section]
\newtheorem{prop}[theorem]{Proposition}

\newtheorem{lemma}[theorem]{Lemma}
\newtheorem{define}[theorem]{Definition}

\newtheorem{example}[theorem]{Example}

\title{
The Shi arrangements and
the Bernoulli polynomials}
\author{Daisuke Suyama
\footnote
{
email: s093014@math.sci.hokudai.ac.jp
}
\,
Hiroaki Terao
\footnote
{
Corresponding author.
Supported by JSPS Grants-in-Aid, Scientific Research
(B) 
No. 21340001.
email: terao@math.sci.hokudai.ac.jp
}
\\
\\
{\footnotesize {\it Department of Mathematics, Hokkaido University, 
Sapporo, Hokkaido 060-0810, Japan.}}
}

\date{}

\begin{document}

\maketitle

\begin{abstract}
The braid arrangement is the
Coxeter arrangement 
of the type $A_\ell$.  
The Shi arrangement is an affine 
arrangement of hyperplanes 
consisting of the hyperplanes 
of the braid arrangement and their 
parallel translations.
In this paper, 
we give an explicit basis
construction 
for the derivation module
of the cone
over the Shi arrangement.  The essential ingredient of our recipe 
is
the Bernoulli polynomials.
\end{abstract} 


\section{Introduction} 

Let $E$ be an $\ell$-dimensional Euclidean space and $\Phi$ be the root system 
of the type $\mathit{A}_{\ell}$. 
Let $\Phi_{+}$ denote the set of positive roots. 
In this paper we explicitly choose $E$ and $\Phi$ as follows:
let
$W = \R^{\ell+1}$
and 
$x_{1}, \dots, x_{\ell+1}$
be an orthonormal basis for the dual space $W^{*}$.
Define 
\begin{align*} 
E &:= \{\sum_{i=1}^{\ell+1} c_{i} x_{i} 
\in W^{*}  \mid \sum_{i=1}^{\ell+1} c_{i} = 0\},\\
\Phi &:= \{x_{i} - x_{j} \in E
\mid 1\leq i\leq \ell+1, 1\leq j\leq \ell+1, 
i\neq j\},
\\
\Phi_{+} &:=
\{
x_{i} - x_{j} \in \Phi \mid i < j 
\}.
\end{align*} 
Let
$\mathcal{A}(\Phi)=\{ \mathit{H}_{\alpha} \mid \alpha \in \Phi_{+} \}$, 
where $\mathit{H}_{\alpha}=\{ v \in \mathit{W} \mid \alpha(v)=0 \}$. 
Then
$\A(\Phi)$ 
is called
a {\bf
braid arrangement}, 
which is undoubtedly
the most-studied arrangement of hyperplanes
in various contexts.
In the study of the Kazhdan-Lusztig
representation theory of the affine Weyl groups,
J.-Y. Shi introduced 
the \textbf{Shi\ arrangements} 
in \cite{shi1} 
as follows:
Let
\[
\mathit{H}_{\alpha,1} = \{ v \in \mathit{W} \mid \alpha(v)=1 \}.
\]
Then the Shi arrangement is given by
\[
\mathcal{A}(\Phi) \cup \{ \mathit{H}_{\alpha,1} \mid \alpha \in \Phi_{+} \} = \bigcup_{\alpha \in \Phi_{+}} \{ \mathit{H}_{\alpha},
\mathit{H}_{\alpha,1} \}.
\]
Embed the $(\ell+1)$-dimensional space $W$ into $V=\R^{\ell+2} 
$ by 
adding a new coordinate $z$ such that $W$ is defined by 
the equation $z=1$ in $V$. 
Then, as in \cite[Definition 1.15]{OT},
we have the {\bf cone}  $\Shi_{\ell} $ over the Shi arrangement.
It is a central arrangement in $V$
defined by
$$Q(\Shi_{\ell} ) = z \prod_{1\leq p<q\leq \ell+1} (x_{p} - x_{q})
\prod_{1\leq p<q\leq \ell+1} (x_{p} - x_{q}-z)=0.$$
Let $\mathit{S}$ be the algebra of polynomial functions on $\mathit{V}$ and let $\mathrm{Der}_{\mathit{V}}$ be the module of derivations 
of $\mathit{S}$ to itself
\[
\mathrm{Der}_{\mathit{V}}=\{ \theta : \mathit{S} \rightarrow \mathit{S} \mid \theta\ \text{is}\ \mathbb{R} \text{-linear and}\ \theta (fg) =f\theta(g)+ g\theta(f)\ \text{for any}\ f,g \in \mathit{S} \}.
\]
The derivation module $\mathit{D}(\mathcal{S}_{\ell})$ is defined by 
\begin{multline*} 
\mathit{D}(\mathcal{S}_{\ell})=\{ \theta \in \mathrm{Der}_{\mathit{V}} \mid 
\theta(\alpha)\ \text{is divisible by}\ \alpha\\
\text{ and }\ 
\theta(\alpha-z)\ \text{is divisible by}\ \alpha-z
\text{ for any}\ \alpha \in \Phi_{+} \}.
\end{multline*} 
In the present paper, for the first time,
we construct an explicit basis for the derivation module 
$\mathit{D}(\mathcal{S}_{\ell})$.
The most important ingredients of our recipe are 
the \textbf{Bernoulli polynomials} $\mathit{B}_{k}(x)$
and their relatives $B_{p, q}(x)$ in Definition \ref{Bpqdefinition}.
The explicit costruction of our basis is in Definition \ref{basisdefinition}.

One of the remarkable properties of the Shi arrangement  is the fact that its number of chambers is equal to
$(\ell+2)^{\ell}$.  
A good number of articles,
including \cite{St96, ER, Hea, Ath, Y04},
study this intriguing property. 
Because of Zaslavsky's chamber counting 
formula \cite{Z75},
the property follows from the formula
\[\pi(\mathcal{S}_{\ell}, t)= (1+t)\left(1+(\ell +1)t\right)^{\ell} \]
for the Poincar\'e polynomial
\cite{OT} of the cone $\mathcal{S}_{\ell}$.    
Ch. Athanasiadis proved that 
$\mathit{D}(\mathcal{S}_{\ell})$ is a free $\mathit{S}$-module 
with exponents
$(0, 1, \ell+1,\ldots ,\ell+1)$ in \cite{Ath}. 
He consequently proved the formula
above
thanks to
the factorization theorem in
\cite{T81} 
which
asserts
that
if the derivation module 
$\mathit{D}(\mathcal{A})$
 is a free $\mathit{S}$-module with a basis 
$\theta_{1},\ldots ,\theta_{\ell}$ 
then the Poincar\'e polynomial
 of 
$\mathcal{A}$ 
is equal to 
$\prod_{j=1}^{\ell}\left(1 +  (\deg \theta_{j})t\right)$.
His proof of the freeness in \cite{Ath} uses the addition-deletion 
theorem \cite{T80}.
Later M.
Yoshinaga extended this result in \cite{Y04}
to the generalized 
Shi and Catalan arrangements and affirmatively settled the Edelman-Reiner 
conjecture
\cite{ER}. However, even in the case of Shi arrangements, 
no basis was constructed so far.

The organization of this paper is as follows:
in Section 2, we will define the polynomials
$B_{p,q}(x) $ which includes the Bernoulli polynomials.
In Section 3,
Theorem \ref{main} proves that the derivations constructed in Definition
\ref{basisdefinition} form
a basis for the derivation module $D(\Shi_{\ell} )$.

\section{The Bernoulli polynomials and $B_{p,q}(x)$}

Let $\mathit{B}_{k} (x)$ 
denote 
the $k$-th
{\bf Bernoulli polynomial}.
Let
$\mathit{B}_{k}(0)=\mathit{B}_{k}$
denote the $k$-th
{\bf Bernoulli number}.
The most important property of the Bernoulli polynomial
in this paper
is the following elementary formula (e.g., \cite{APO76}):

\begin{theorem}
\label{bernoullipolynomialdifference}
\[
B_{k} (x+1) - B_{k}(x) = k x^{k-1}.  
\]
\end{theorem}

\begin{define}
\label{Bpqdefinition} 
For $(p,q) \in (\mathbb{Z}_{\geq 0})^2$, consider a polynomial
$\mathit{B}_{p,q} (x)$ in $x$ satisfying the following two conditions:
\begin{enumerate}
\item $\mathit{B}_{p,q} (x+1)- \mathit{B}_{p,q} (x)=(x+1)^{p} x^{q}$,
\item $\mathit{B}_{p,q} (0)=0$.
\end{enumerate}
It is easy to see that $\mathit{B}_{p,q} (x)$
 is uniquely determined by these two conditions.
\end{define}

\begin{example}
\label{examplebpq} 
(1)
When $(p,q)=(0,q)$, we have
$$\mathit{B}_{0,q} (x)=\frac{1}{q+1}
 \left\{
\mathit{B}_{q+1}(x)-\mathit{B}_{q+1}\right\}$$ 
because of 
Theorem \ref{bernoullipolynomialdifference}. 

(2)
When $(p,q)=(p,0)$, we obtain
$$\mathit{B}_{p,0}(x)
=
\frac{(-1)^{p+1}}{p+1}
\left\{
\mathit{B}_{p+1}(-x)
-
\mathit{B}_{p+1}
\right\}
=
(-1)^{p+1} B_{0,p}(-x) 
$$ 
because
\begin{align*}
&~~~~
(-1)^{p+1} B_{0,p}(-x-1)
-
(-1)^{p+1} B_{0,p}(-x)\\
&=
(-1)^{p} 
\left\{
B_{0,p}(-x)
-
B_{0,p}(-x-1)
\right\}
= (-1)^{p} (-x-1)^{p} 
= (x+1)^{p}. 
\end{align*}

(3)
For a general $(p,q) \in (\mathbb{Z}_{\geq 0})^{2}$, 
it easily follows from
Theorem \ref{bernoullipolynomialdifference}
that the polynomial has an expression in terms of the
Bernoulli polynomials as
\begin{equation*}
\mathit{B}_{p,q}(x)=
\sum_{i=0}^{p} 
\frac{1}{q+i+1} 
\binom{p}{i}
\left\{
\mathit{B}_{q+i+1}(x)-\mathit{B}_{q+i+1}
\right\}
=
\sum_{i=0}^{p} 
\binom{p}{i}
B_{0,q+i} (x).
\end{equation*}
For example, $B_{1,1}(x) = B_{0,1}(x) + B_{0,2}(x) =
\frac{1}{3}   (x^{3} - x) $. 
\end{example}

Note that the polynomial $B_{p,q} (x)$ 
is a polynomial of degree $p+q+1$. 
The homogenization $\overline{B}_{p,q} (x, z) $ 
of $B_{p, q} (x) $ is defined by
\[
\overline{B}_{p,q} (x, z) := 
z^{p+q+1} B_{p,q} \left(\frac{x}{z}\right).
\]

\section{A basis construction} 

Let $1\leq j\leq \ell$. Define 
$$\mathit{I}_{1}= \{ x_{1},x_{2},\ldots ,x_{j-1} \},
\,\,\,
\mathit{I}_{2}= \{ x_{j+2},x_{j+3},\ldots ,x_{\ell +1} 
\}.$$ 
Let $\sigma_{k}^{(s)}$ denote the elementary 
symmetric function 
in the variables in
$\mathit{I}_{s}$ of degree $k$ 
$(s=1,2\ ,\ k \in \mathbb{Z}_{\geq 0})$. 
Recall the homogeneous
polynomials
$\overline{B}_{p,q} (x, z)$ 
of degree $p+q+1$ defined at the end of the previous section.

\begin{define}
\label{basisdefinition} 
Let
$\partial_{i} \,\,
(1\leq i\leq \ell+1)
$ and $\partial_{z} $ denote 
$\partial/\partial x_{i} $ 
and
$\partial/\partial z $
respectively. 
Define homogeneous derivations 
$$
\eta_{1} := \sum_{i=1}^{\ell+1} \partial_{i}
\in D(\Shi_{\ell}),
\,\,\,\,
\eta_{2} :=  z \partial_{z} +\sum_{i=1}^{\ell+1} x_{i} \partial_{i}
\in D(\Shi_{\ell}),
$$ 
and
\begin{equation*}
\varphi_{j}
:=
(x_{j}-x_{j+1}-z) 
\sum_{i=1}^{\ell +1} \sum_{\substack{
0\leq k_{1} \leq j-1 \\ 0\leq k_{2} \leq \ell-j}} 
(-1)^{k_{1} +k_{2} } 
\sigma_{j-1-k_{1}}^{(1)} 
\sigma_{\ell-j-k_{2}}^{(2)} 
\,
\overline{B}_{k_{1}, k_{2}}
(x_{i}, z) \partial_{i}
\end{equation*}
for
$1\leq j \leq \ell$.
\end{define}

We will prove
that  the derivations
$
\eta_{1}, 
\eta_{2},$
and
$ 
\varphi_{1}, \dots , \varphi_{\ell}$
form a basis for $D(\Shi_{\ell})$. 
First we will verify the following Proposition:

\begin{prop}
\label{Proposition6} 
The derivations 
$\varphi_{j}
\,\,
(1\leq j\leq \ell)
$ 
belong to the module
$D(\Shi_{\ell} )$. 
\end{prop}

\noindent
{\textbf{proof.}}
We first have
\begin{align*}
\varphi_{j}(x_{p}-x_{q}) = &
(x_{j}-x_{j+1}-z) \\
&~\sum_{\substack{0\leq k_{1}\leq j-1 \\ 0\leq k_{2}\leq \ell-j}} 
(-1)^{k_{1} +k_{2}} 
\sigma_{j-1-k_{1}}^{(1)} 
\sigma_{\ell-j-k_{2}}^{(2)} 
\left\{
\overline{B}_{k_{1},k_{2}}(x_{p}, z) 
- 
\overline{B}_{k_{1}, k_{2}}(x_{q}, z) 
\right\}. 
\end{align*}
Since the right hand side equals zero if we set $x_{p}=x_{q}$, 
we may conclude that $\varphi_{j} (x_{p}-x_{q})$
  is divisible by 
$x_{p} - x_{q} $ for all pairs $(p,q)$ with
$1\leq p < q\leq \ell+1$. 

The congruent notation $\equiv$ in the following 
calculation is modulo the ideal $(x_{p}-x_{q}-z)$: 
\begin{align*}  
&~~~\varphi_{j} (x_{p}-x_{q} -z)\\
&\equiv (x_{j}-x_{j+1}-z)\\
&~~~\sum_{\substack{0\leq k_{1}\leq j-1 \\ 0\leq k_{2}\leq \ell-j}} 
(-1)^{k_{1} +k_{2}} 
\sigma_{j-1-k_{1}}^{(1)} 
\sigma_{\ell-j-k_{2}}^{(2)} 
\left\{
\overline{B}_{k_{1}, k_{2}}(x_{p}, x_{p} - x_{q}) 
- 
\overline{B}_{k_{1}, k_{2}}(x_{q}, x_{p} - x_{q}) 
\right\}\\
&=(x_{j}-x_{j+1}-z) 
\sum_{\substack{
 0\leq k_{1}\leq j-1\\
0\leq k_{2}\leq \ell-j}} 
(-1)^{k_{1} +k_{2} } 
\sigma_{j-1-k_{1}}^{(1)} \sigma_{\ell-j-k_{2}}^{(2)} \\
& ~~~~~~~~~~~~~~~~~~~~~~~~~~~~~~~~~
(x_{p}-x_{q})^{k_{1} + k_{2}+1}
\{ \mathit{B}_{k_{1}, k_{2}}
(\frac{x_{p}}{x_{p}-x_{q}}) - \mathit{B}_{k_{1}, k_{2}}
(\frac{x_{q}}{x_{p}-x_{q}}) \} \\
&=(x_{j}-x_{j+1}-z) 
 \\
&~~~~~~~~~
\sum_{\substack{0\leq k_{1}
\leq j-1 \\ 0\leq k_{2}\leq \ell-j}} 
(-1)^{k_{1} +k_{2} } 
\sigma_{j-1-k_{1}}^{(1)} \sigma_{\ell-j-k_{2}}^{(2)} 
(x_{p}-x_{q})^{k_{1} + k_{2}+1 } (\frac{x_{p}}{x_{p}-x_{q}})^{k_{1}}
 (\frac{x_{q}}{x_{p}-x_{q}})^{k_{2}} \\
&= (x_{j}-x_{j+1}-z) (x_{p}-x_{q}) 
\sum_{\substack{
0\leq k_{1}\leq j-1 \\ 0\leq k_{2}\leq \ell-j}} 
(-1)^{k_{1} +k_{2} } 
\sigma_{j-1-k_{1}}^{(1)} \sigma_{\ell-j-k_{2}}^{(2)} 
x_{p}^{k_{1}} x_{q}^{k_{2}} \\
&= (x_{j}-x_{j+1}-z) (x_{p}-x_{q}) \sum_{k_{1}=0}^{j-1} 
\sigma_{j-1-k_{1}}^{(1)} 
(-x_{p})^{k_{1}} 
\sum_{k_{2}=0}^{\ell -j} 
\sigma_{\ell-j-k_{2}}^{(2)} (-x_{q})^{k_{2}}\\ 
&= (x_{j}-x_{j+1}-z) (x_{p}-x_{q}) 
\prod_{s=1}^{j-1}(x_{s}-x_{p}) \prod_{s=j+2}^{\ell +1}(x_{s}-x_{q}) 
\equiv 
0 
\end{align*}
for all pairs $(p, q)$ with $1\leq p < q \leq \ell +1$. 
$\square$

\begin{lemma}
\label{detSigma}
Suppose $\ell\geq 1$. 
Let $N$ be the $\ell\times\ell$-matrix whose
$(i, j)$-entry is equal to the elementary symmetric function of degree 
$\ell-i$ 
in the variables $x_{1}, \dots, x_{j-1}, x_{j+2}, \dots, x_{\ell+1   } $.
Then 
\[
\det N
=
(-1)^{\ell(\ell-1)/2} 
\prod_{\substack{1\leq p<q\leq\ell\\ q-p>1}} (x_{p} - x_{q} ).
\]
 \end{lemma}

\noindent
{\textbf{proof.}}
Note that we have the equality
\begin{align*} 
&~~~
\left[
1 \,\,\, -x_{p} \,\,\,\, (-x_{p})^{2} \dots  (-x_{p})^{\ell-2}
\,\,\,\, (-x_{p})^{\ell-1}
\right]
N\\
&=
\left[
\prod_{\substack{1\leq s\leq \ell+1\\
s\not\in \{1,2\}}}
(x_{s} -  x_{p})
\, 
\prod_{\substack{1\leq s\leq \ell+1\\
s\not\in \{2,3\}}}
(x_{s} -  x_{p})
\, 
\dots
\,
\prod_{\substack{1\leq s\leq \ell+1\\
s\not\in \{\ell-1,\ell\}}}
(x_{s} -  x_{p})
\,
\prod_{\substack{1\leq s\leq \ell+1\\
s\not\in \{\ell,\ell+1\}}}
(x_{s} -  x_{p})
\right]
\end{align*} 
for any $1\leq p\leq \ell$. 
Suppose that
$$1\leq p < q \leq \ell+1,
\,\,\,
q-p>1.$$
Set $x_{p} = x_{q} $ in $N$, and we get
$N_{pq} $.
Then we may conclude that
\begin{align*} 
&~~~
\left[
1 \,\,\, -x_{p} \,\,\,\, (-x_{p})^{2} \dots  (-x_{p})^{\ell-2}
\,\,\,\, (-x_{p})^{\ell-1}
\right]
N_{pq}
={\mathbf 0}.
\end{align*} 
This implies that
$\det N_{pq} =0$
and that
$
\det N
$
is divisible by 
$x_{p} - x_{q}.$ 
Since
\[
\deg 
(\det
N)
=
\ell(\ell-1)/2
=
\deg
\prod_{\substack{1\leq p<q\leq\ell+1\\ q-p>1}} 
(x_{p} - x_{q}),
\]
there exists a constant $C$ such that
\[
\det N
=
C\,
(-1)^{\ell(\ell-1)/2} 
\prod_{\substack{1\leq p<q\leq\ell+1\\ q-p>1}} 
(x_{p} - x_{q})
=
C
\prod_{\substack{1\leq p<q\leq\ell+1\\ q-p>1}} 
(x_{q} - x_{p}).
\]
By comparing the coefficients of 
$
x_{3} x_{4}^{2}   \dots x_{\ell}^{\ell-2} x_{\ell+1}^{\ell-1}     
$ 
on both sides, we obtain $C=1$. 
$\square$

\begin{prop}
\label{Proposition8} 
The derivations 
$\eta_{1}, $
$\eta_{2}, $
$\varphi_{1} , \dots , \varphi_{\ell}$ 
are linearly independent over $S$.  
\end{prop}

\noindent
{\textbf{proof.}}
Set $z=0$ in $\varphi_{j} $ and we get
$\phi_{j} $ as follows:
\begin{align*} 
\phi_{j} &:= \varphi_{j} |_{z=0}
=
(x_{j}-x_{j+1})
\sum_{i=1}^{\ell +1} \sum_{\substack{
0\leq k_{1} \leq j-1 \\ 0\leq k_{2} \leq \ell-j}} 
\frac{(-1)^{k_{1} +k_{2} }}{k_{1} +k_{2} +1} 
\sigma_{j-1-k_{1}}^{(1)} 
\sigma_{\ell-j-k_{2}}^{(2)} 
\,
x_{i}^{k_{1} +k_{2} +1}  \partial_{i}\\
&= 
(x_{j}-x_{j+1})
\sum_{k=1}^{\ell}  
\frac{(-1)^{k-1}}{k}
\left(
\sum_{\substack{k_{1} +k_{2} +1 =k\\
0\leq k_{1} \leq j-1 \\ 0\leq k_{2} \leq \ell-j}} 
\sigma_{j-1-k_{1}}^{(1)} 
\sigma_{\ell-j-k_{2}}^{(2)} 
\right)
\sum_{i=1}^{\ell +1} 
x_{i}^{k}   
\partial_{i}\\
&=
(x_{j}-x_{j+1})
\sum_{k=1}^{\ell}  
\frac{(-1)^{k-1}}{k}\,
\sigma_{\ell-k}
(x_{1}, \dots, x_{j-1}, x_{j+2}, \dots, x_{\ell+1})
\,
\sum_{i=1}^{\ell +1} 
x_{i}^{k}   
\partial_{i}.
\end{align*} 
Here
$
\sigma_{\ell-i}
(x_{1}, \dots, x_{j-1}, x_{j+2}, \dots, x_{\ell+1})
$ 
stands for
the elementary symmetric function
of degree $\ell-i$ in the variables 
$x_{1}, \dots, x_{j-1}, x_{j+2}, \dots, x_{\ell+1}.
$ 
This is equal to 
the $(i, j)$-entry 
$N_{ij} $ of the matrix
$
N
$  
in Lemma \ref{detSigma}.
Thus we have
\begin{equation} 
\phi_{j} (x_{i}) =
(x_{j}-x_{j+1})
\sum_{k=1}^{\ell}  
\frac{(-1)^{k-1}}{k}\,
x_{i}^{k}   
N_{kj}. 
\label{phijxi} 
\end{equation} 
Define
 two $(\ell+1)\times(\ell+1)$-diagonal
 matrices
$D_{1} $ and
$D_{2} $ 
by
\begin{align*} 
D_{1} 
&:= 
\left[1\right]\oplus\left[1\right]\oplus\left[(-1)^{1}/2
\right]\oplus \left[(-1)^{2}/3\right] 
\oplus\dots\oplus [(-1)^{\ell-1}/\ell],
\\
D_{2} 
&:=
[1]\oplus[x_{1} -x_{2}]\oplus[x_{2} - x_{3}]
\oplus 
\dots
\oplus
[x_{\ell}-x_{\ell+1}],
\end{align*} 
where $\oplus$ stands for the direct sum of matrices.
Also define
 two $(\ell+1)\times(\ell+1)$-matrices
$\tilde{N} $ and
$M$ 
by
$$
\tilde{N} := [1]\oplus N,
\,\,\,\,\,
M:= 
\left[
x_{i}^{j-1}  
\right]_{
1\leq i\leq\ell+1,
1\leq j\leq\ell+1
} .
$$
From (\ref{phijxi})
 we obtain
\[
P
 :=
\begin{bmatrix}
1&\phi_{1}(x_{1}) &\dots&\phi_{\ell}(x_{1})\\  
1&\phi_{1} (x_{2}) &\dots&\phi_{\ell}(x_{2})\\  
1&\phi_{1} (x_{3}) &\dots&\phi_{\ell}(x_{3})\\  
.&  .  & \dots& .\\  
.&  .  & \dots& .\\  
1&\phi_{1} (x_{\ell+1}) &\dots&\phi_{\ell} (x_{\ell+1})
\end{bmatrix}
=
M
D_{1} 
\tilde{N} 
D_{2}.
\]
Thus, by applying the Vandermonde determinant
formula
and Lemma \ref{detSigma}, we deduce
\begin{align*} 
\det P
&=
(\det M)
(\det D_{1}) 
(\det \tilde{N}) 
(\det D_{2})
\\
&=
\left(
\prod_{1\leq p < q\leq \ell+1} (x_{q}-x_{p})
\right)
\left(
\frac{(-1)^{\ell(\ell-1)/2}}{\ell !} 
\right)
\left(
\det {N} 
\right)
\prod_{1\leq p \leq \ell} (x_{p}-x_{p+1})
\\  
&=
\left(
\frac{(-1)^{\ell(\ell+1)/2}}{\ell !} 
\right)
\prod_{1\leq p < q \leq \ell+1} (x_{p}-x_{q})^{2}
\neq 0.  
\end{align*} 
Thus 
$\eta_{1}, \phi_{1}, \dots, \phi_{\ell}$
are linearly independent.  
This implies that
$\eta_{1}, \varphi_{1}, \dots, \varphi_{\ell}$
are linearly independent. 
Since $\eta_{2}(z) = z $ 
and
$\eta_{1}(z)=\varphi_{1}(z)= \dots =\varphi_{\ell}(z)=0$,
we finally conclude that
$\eta_{1}, \eta_{2}, \varphi_{1}, \dots, \varphi_{\ell}$
are linearly independent. 
$\square$

\bigskip

\noindent
{\em Remark.} 
The derivations 
$\phi_{1}, \dots, \phi_{\ell}  $ are a basis
for the derivation module of the double 
Coxeter arrangement of the
type $A_{\ell}$ studied in \cite{ST98}
(cf. \cite{T02}).

\begin{theorem}
\label{main} 
The derivations $\eta_{1} , \eta_{2}, \varphi_{1}, \dots, \varphi_{\ell}   $ 
form a basis for $D(\Shi_{\ell})$. 
\label{main} 
\end{theorem}

\noindent
{\textbf{proof.}}
We may apply
Saito's criterion 
\cite{S80} \cite[Theorem 4.23]{OT} 
thanks to Propositions \ref{Proposition6}
and
 \ref{Proposition8} 
because 
\[
\deg \eta_{1} + \deg \eta_{2} +\sum_{j=1}^{\ell}  \deg \varphi_{j} 
=
1+\ell(\ell+1)=|\Shi_{\ell} |.
\,\,\,
\square
\]

\noindent
{\textit{Remark.}}
The Bernoulli polynomials explicitly appear
in the first derivation $\varphi_{1} $ and the last one
$\varphi_{\ell} $ because of Example \ref{examplebpq} (1) and (2):
\begin{align*}
\varphi_{1}
&=
(x_{1}-x_{2}-z) 
\sum_{i=1}^{\ell +1} \sum_{k_{2} =0}^{\ell-1} 
(-1)^{k_{2} } 
\sigma_{\ell-1-k_{2} }^{(2)} 
\,
\overline{B}_{0,k_{2} }
(x_{i}, z) \partial_{i}\\
&=
(x_{1}-x_{2}-z) 
\sum_{i=1}^{\ell +1} 
\sum_{k=1}^{\ell} 
\frac{(-1)^{k-1}}{k}  
\sigma_{\ell-k}^{(2)} 
z^{k} 
\left(
B_{k}(x_{i}/z)
-
B_{k}
\right) 
 \partial_{i},
\end{align*}
and
\begin{align*}
\varphi_{\ell}
&=
(x_{\ell}-x_{\ell+1}-z) 
\sum_{i=1}^{\ell +1} \sum_{k_{1} =0}^{\ell-1} 
(-1)^{k_{1} } 
\sigma_{\ell-1-k_{1} }^{(1)} 
\,
\overline{B}_{k_{1} ,0}
(x_{i}, z) \partial_{i}\\
&=
(x_{\ell}-x_{\ell+1}-z) 
\sum_{i=1}^{\ell +1} \sum_{k=1}^{\ell} 
\frac{(-1)^{k-1}}{k}  
\sigma_{\ell-k}^{(1)} 
\,
(-z)^{k} 
\left(
B_{k}(-x_{i}/z) 
-
B_{k}
\right)
 \partial_{i}.
\end{align*}
Here 
$\sigma^{(1)}_{d}$
 and
$\sigma^{(2)}_{d}$
are the elementary
symmetric functions
of degree $d$ in the variables 
$x_{1}, \dots, x_{\ell-1}  $ 
and
$x_{3}, \dots, x_{\ell+1}  $ 
respectively.

\begin{example}
For $A_{3}$, we have
\begin{align*}
\eta_{1} &= \partial_{1} +\partial_{2} +\partial_{3} +\partial_{4}, \,\,\,\,\,\,\,\,
\eta_{2} = x_{1} \partial_{1} +x_{2} \partial_{2} +x_{3} \partial_{3} +x_{4} \partial_{4} +z \partial_{z}, \\
\varphi_{1}&=
%\frac{1}{6}
x_{1}(x_{1}-x_{2}-z)
%
%(
%6x_{3}x_{4}
%-3x_{3}x_{1}
%+3x_{3}z
%-3x_{4}x_{1}
%+3x_{4}z
%+2x_{1}^{2}
%-3x_{1}z
%+z^{2}
%)
\left\{
x_{3}x_{4}
-
\frac{1}{2}
(x_{3}+x_{4})
(x_{1}
-z)
+
\frac{1}{3}
\left(x_{1}^{2}
-
\frac{3}{2}x_{1}z
+
\frac{1}{2}
z^{2}\right)
\right\}
\partial_{1} \\
&\, +
%\frac{1}{6}
x_{2}
(x_{1}-x_{2}-z)
%
%(
%6x_{3}x_{4}
%-3x_{2}x_{3}
%-3x_{2}x_{4}
%+3x_{3}z
%+3x_{4}z
%+2x_{2}^{2}
%-3x_{2}z
%+z^{2}
%)
\left\{
x_{3}x_{4}
-
\frac{1}{2}
(x_{3}+x_{4})
(x_{2}
-z)
+
\frac{1}{3}
\left(x_{2}^{2}
-
\frac{3}{2}x_{2}z
+
\frac{1}{2}
z^{2}\right)
\right\}
\partial_{2} \\
&\ \ -\frac{1}{6}x_{3}(x_{1}-x_{2}-z)(x_{3}+z)(x_{3}-3x_{4}-z)\partial_{3} \\
&\ \ -\frac{1}{6}x_{4}(x_{1}-x_{2}-z)(x_{4}+z)(x_{4}-3x_{3}-z)\partial_{4},
\\
%
%\end{align*}
%\begin{align*}
\varphi_{2}&=-\frac{1}{6}x_{1}(x_{2}-x_{3}-z)(x_{1}-z)(x_{1}-3x_{4}-2z)
\partial_{1} \\
&\ \ +
%\frac{1}{6}
x_{2}(x_{2}-x_{3}-z)
\left\{
x_{1}x_{4}
-
\frac{1}{2}x_{1}
(x_{2}-z)
-
\frac{1}{2}x_{4}
(x_{2}+z)
+
\frac{1}{3}
(x_{2}^{2}
-z^{2}
)\right\}
\partial_{2} \\
&\ \ +
%\frac{1}{6}
x_{3}(x_{2}-x_{3}-z)
\left\{
x_{1}x_{4}
-
\frac{1}{2}x_{1}
(x_{3}-z)
-
\frac{1}{2}x_{4}
(x_{3}+z)
+
\frac{1}{3}
(x_{3}^{2}
-z^{2}
)\right\}
%
%(
%6x_{1}x_{4}
%-3x_{1}x_{3}
%+3x_{1}z
%-3x_{3}x_{4}
%-3x_{4}z
%+2x_{3}^{2}
%-2z^{3}
%)
\partial_{3} \\
&\ \ +\frac{1}{6}x_{4}(x_{2}-x_{3}-z)(x_{4}+z)(3x_{1}-x_{4}-2z)\partial_{4},
\\
\varphi_{3}&
=-\frac{1}{6}x_{1}(x_{3}-x_{4}-z)(x_{1}-z)(x_{1}-3x_{2}+z)\partial_{1} \\
&\ \ \, -\frac{1}{6}x_{2}(x_{3}-x_{4}-z)(x_{2}-z)(x_{2}-3x_{1}+z)\partial_{2} \\
&\, +
%\frac{1}{6}
x_{3}(x_{3}-x_{4}-z)
\left\{
x_{1}x_{2}
-
\frac{1}{2}
(x_{1}+x_{2})
(x_{3}
+z)
+
\frac{1}{3}
\left(x_{3}^{2}
+
\frac{3}{2}x_{3}z
+\frac{1}{2}
z^{2}\right)
\right\}
%
%(
%6x_{1}x_{2}
%-3x_{1}x_{3}
%-3x_{1}z
%-3x_{2}x_{3}
%-3x_{2}z
%+2x_{3}^{2}
%+3x_{3} z 
%+z^{2}
%)
\partial_{3} \\
&\,
+
%\frac{1}{6}
x_{4}(x_{3}-x_{4}-z)
\left\{
x_{1}x_{2}
-
\frac{1}{2}
(x_{1}+x_{2})
(x_{4}+z)
+
\frac{1}{3}
\left(x_{4}^{2}
+
\frac{3}{2}x_{4}z+
\frac{1}{2}
z^{2}\right)
\right\}
%
%(
%6x_{1}x_{2}
%-3x_{1}x_{4}
%-3x_{1}z
%-3x_{2}x_{4}
%-3x_{2}z
%+2x_{4}^{2}
%+3x_{4} z 
%+z^{2}
%)
\partial_{4}.
\end{align*}
\end{example}

\bigskip

\noindent
{\textit{Problem.}}
Construct a basis for the derivation module
$D(\A)$ when $\A$ is the cone over the generalized 
Shi arrangement of one
of the other types: $B_{\ell}, C_{\ell}, 
D_{\ell}, E_{6}, E_{7}, E_{8},
F_{4}$ and $G_{2}.$ 
It has been known that the derivation modules 
are free modules for these cases by M. Yoshinaga
\cite{Y04}.  It seems interesting to learn what 
kind of polynomials plays the role of
the Bernoulli polynomials in the case of
the type $A_{\ell}$.

 \vspace{5mm}

\end{document}